%% file: Multichannel-LRD.tex
\documentclass[11pt]{article}
\usepackage{graphicx}
\usepackage{amsmath}
\usepackage{amsfonts}
\usepackage{epsfig}

\newcommand{\be}{\begin{equation}}
\newcommand{\ee}{\end{equation}}
\newcommand{\bes}{\begin{equation*}}
\newcommand{\ees}{\end{equation*}}
\newcommand{\beqn}{\begin{eqnarray}}
\newcommand{\eeqn}{\end{eqnarray}}
\newcommand{\beqns}{\begin{eqnarray*}}
\newcommand{\eeqns}{\end{eqnarray*}}

\newcommand{\fr}[1]{(\ref{#1})}

\newcommand{\lkr}{\left(}
\newcommand{\rkr}{\right)}
\newcommand{\lkv}{\left[}
\newcommand{\rkv}{\right]}
\newcommand{\lfi}{\left\{}
\newcommand{\rfi}{\right\}}

\newcommand{\EE}{\ensuremath{{\mathbb E}}}
\newcommand{\II}{\ensuremath{{\mathbb I}}}
\newcommand{\PP}{\ensuremath{{\mathbb P}}}

\newcommand{\RR}{\ensuremath{{\mathbb R}}}
\newcommand{\CC}{\ensuremath{{\mathbb C}}}

\newcommand{\ints}{\ensuremath{{\mathbb Z}}}

\newcommand{\sumku}{\sum_{k \in \Ujr}}
\newcommand{\sumr}{\sum_{r \in \Aj}}

\newcommand{\card}{\mbox{card}}
\newcommand{\Var}{\mbox{Var}}
\newcommand{\Cov}{\mbox{Cov}}

\newcommand{\psijk}{\psi_{jk}}
\newcommand{\ajk}{a_{j_0k}}
\newcommand{\bjk}{b_{jk}}
\newcommand{\hajk}{\widehat{a}_{j_0k}}
\newcommand{\hbjk}{\widehat{b}_{jk}}

\newcommand{\fm}{f_m}
\newcommand{\hfm}{\widehat{f}_m}
\newcommand{\psimjk}{\psi_{mjk}}
\newcommand{\phmjk}{\ph_{mj_0k}}

\newcommand{\Aj}{A_j}
\newcommand{\Ujr}{U_{jr}}
\newcommand{\Bjr}{B_{jr}}
\newcommand{\hBjr}{\widehat{B}_{jr}}

\newcommand {\bSigma}{\mbox{\mathversion{bold}$\Sigma$}}

\newcommand {\bxi}{\mbox{\mathversion{bold}$\xi$}}
\newcommand {\bomega}{\mbox{\mathversion{bold}$\omega$}}
\newcommand {\bA}{\mathbf{A}}

\newcommand {\bV}{\mathbf{V}}

\newcommand {\bh}{\mathbf{h}}

\newcommand{\eps}{\varepsilon}
\newcommand{\ph}{\varphi}

\newcommand{\te}{\theta}
\newcommand{\lam}{\lambda}

\newcommand{\jo}{{j_0}}

\newcommand{\ns}{{n^*}}

\newtheorem{theorem}{Theorem}
\newtheorem{lemma}{Lemma}
\newtheorem{corollary}{Corollary}

\newtheorem{remark}{Remark}
\newtheorem{example}{Example}

\newcommand{\dus}{d^{*}}

\begin{document}

\title{\bf {Multichannel Deconvolution with Long-Range Dependence: A Minimax Study}}

\author{{\em Rida Benhaddou$^1$, Rafal Kulik$^2$, Marianna Pensky$^1$ and Theofanis Sapatinas$^3$ \vspace{0.2cm}} \\
$^1$ Department of Mathematics, University of Central Florida, USA\\
$^2$ Department of Mathematics and Statistics, University of Ottawa, Canada\\
$^3$ Department of Mathematics and Statistics, University of Cyprus, Cyprus}
\date{}

\bibliographystyle{plain}
\maketitle

\begin{abstract}
We consider the problem of estimating the unknown response function in the  multichannel deconvolution model with  long-range dependent Gaussian errors.
We do not limit our consideration to a specific type of long-range dependence rather we assume that the errors should satisfy a general assumption in terms of the smallest and larger eigenvalues of their covariance matrices. We derive minimax lower bounds for the quadratic risk in the proposed
multichannel deconvolution model when the response function is assumed to belong to a Besov ball and the blurring function is assumed to possess some
smoothness properties, including both regular-smooth and
super-smooth convolutions. Furthermore, we propose an adaptive
wavelet estimator of the response function that is asymptotically optimal (in
the minimax sense), or near-optimal within a logarithmic factor, in
a wide range of Besov balls. It is shown that the optimal convergence rates depend on the balance between the smoothness parameter of the response function, the kernel parameters of the blurring function, the long memory parameters of the errors, and how the total number of observations is distributed among the total number of channels. Some examples
of inverse problems in mathematical physics where one needs to recover initial
or boundary conditions on the basis of observations from a noisy solution of a partial differential equation are used to illustrate the application of the theory we developed. The optimal convergence rates and the adaptive estimators we consider extend the ones studied by Pensky and Sapatinas (2009, 2010) for independent and identically distributed Gaussian errors to the case of long-range dependent Gaussian errors. \\

\noindent
{\bf AMS (2000) Subject Classifications:} 62G05 (primary), 62G08, 35J05, 35K05, 35L05 (secondary)\\

\noindent
{\bf Keywords and Phrases:} adaptivity, Besov spaces, block thresholding,
deconvolution, Fourier analysis, functional data, long-range dependence, Meyer wavelets, minimax estimators, multichannel
deconvolution, partial differential equations, stationary sequences, wavelet
analysis.

\end{abstract}

\section{Introduction}
\label{sec:introduction}
\setcounter{equation}{0}

We consider the estimation problem of the unknown response function $f(\cdot)\in L^2(T)$ from observations $y(u_l, t_i)$ driven by
\be \label{sampl}
y(u_l, t_i) = \int_{T}  g(u_l, t_i -x) f( x) dx  +  \xi_{li},\quad  l=1,2,\ldots,M,\ i=1,2,\ldots,N,
\ee
where $u_l \in U=[a, b]$, $0<a \leq b < \infty$, $T= [0, 1]$, $t_i = i/N$, and the errors $\xi_{li}$ are
Gaussian random variables, independent for different $l$'s, but dependent for different $i$'s.\\

Denote the total number of observations $n=NM$ and assume, without loss of generality, that $N=2^J$
for some integer $J >0$. For each $l=1,2,\ldots, M$, let $\bxi^{(l)}$ be a Gaussian vector with components $ \xi_{li}$,
$ i=1,2,\ldots, N$, and let $\bSigma^{(l)} := \Cov (\bxi^{(l)}) := \EE [ \bxi^{(l)} (\bxi^{(l)})^T ]$ be its covariance matrix.\\

{\bf Assumption A1}:\ \ For each $l=1,2,\ldots, M$, $\bSigma^{(l)}$ satisfies the following condition:  
there exist constants $K_1$ and $K_2$ ($0<K_1 \leq K_2 <\infty$), independent
of $l$ and $N$, such that, for each $l=1,2,\ldots,M$,
\be \label{eq:assump1}
K_1 N^{2d_l} \leq \lambda_{\min} (\bSigma^{(l)}) \leq  \lambda_{\max} (\bSigma^{(l)}) \leq K_2 N^{2d_l}, \quad 0 \leq d_l < 1/2,
\ee
where $\lambda_{\min}$ and $\lambda_{\max}$ are the smallest and the largest eigenvalues of (the Toeplitz matrix) $\bSigma^{(l)}$.
(Here, and in what follows, ``$^T$'' denotes the transpose of a vector or a matrix.)
\\

Assumption A1 is valid when, for each $l=1,2,\ldots, M$, $\bxi^{(l)}$ is a second-order stationary Gaussian sequence with spectral density satisfying
certain assumptions. We shall elaborate on this issue in Section \ref{sec:long-range}.
Note that, in the case of independent errors, for each $l=1,2,\ldots,M$, $\bSigma^{(l)}$ is proportional to the identity matrix and that $d_l =0$.
In this case, the multichannel deconvolution model  \fr{sampl} reduces to the one with independent and identically distributed Gaussian errors.
In a view of \fr{sampl}, the limit situation $d_l =0$, $l=1,2,\ldots,M$,  can be thought of as the standard multichannel deconvolution model
described in Pensky and Sapatinas (2009, 2010). \\


Model \fr{sampl} can also be thought of as the discrete version of a model referred to as the functional deconvolution model 
by Pensky and Sapatinas (2009, 2010).   
The functional deconvolution model has a multitude of applications.
In particular, it can be used in a number of inverse problems in mathematical physics where one needs to recover initial
or boundary conditions on the basis of observations from a noisy solution of a partial differential equation. For instance,
the problem of recovering the initial condition for parabolic equations based on observations in a fixed-time trip was
first investigated in Lattes and Lions~(1967), and the problem of recovering the boundary condition for elliptic equations
based on observations in an interval domain was studied in Golubev and Khasminskii~(1999) and Golubev~(2004).\\

In the case when $a=b$, the functional deconvolution model reduces to the standard deconvolution model.  
This model has been the subject of a great array of research papers since late 1980s,
but the most significant contribution was that of Donoho~(1995) who was the first to device a wavelet
solution to the problem. This has attracted the attention of a good deal of researchers, see, e.g.,
Abramovich and Silverman~(1998), Kalifa and Mallat~(2003), Donoho and Raimondo~(2004), Johnstone and Raimondo~(2004),
Johnstone, Kerkyacharian, Picard and Raimondo~(2004), Kerkyacharian, Picard and Raimondo~(2007). 
(For related results on the density deconvolution problem, we refer
to, e.g., Pensky and Vidakovic (1999), Walter and Shen (1999), Fan and
Koo (2002).)\\

In the multichannel deconvolution model studied by Pensky and Sapatinas~(2009, 2010), as well 
as in the very current extension of their results to derivative estimation by Navarro {\it et al.} (2013), 
it is assumed that errors are independent and identically distributed Gaussian random variables. 
However, empirical evidence has shown  that even at large lags, the correlation structure in the 
errors can decay at a hyperbolic rate, rather than an exponential rate.
To account for this, a great deal of papers on long-range dependence (LRD) have been developed. The study of LRD (also called long memory)  
has a number of  applications, as it can be reflected by the very large number of articles having LRD or long memory
in their titles, in areas such as climate study, DNA sequencing, econometrics, finance, hydrology, internet modeling, signal
and image processing, physics and even linguistics. Other applications can be found in Beran~(1992, 1994), Beran {\it et al.} (2013)
and Doukhan {\em et al}.~(2003). \\

Although quite a few LRD models have been considered in the regression estimation framework,
very little has been done in the standard deconvolution model. The density deconvolution set up has also witnessed
some shift towards analyzing the problem for dependent processes. The argument behind that was that a number of
statistical models, such as non-linear GARCH and continuous-time stochastic volatility models, can be
looked at as density deconvolution models if we apply a simple logarithmic transformation, and thus there is need
to account for dependence in the data. This started by Van Zanten {\em et al}.~(2008) who investigated wavelet based
density deconvolution studied by Pensky and Vidakovic~(1999) with a relaxation to weakly dependent processes.
Comte {\em et al}.~(2008) analyzed another adaptive estimator that was proposed earlier but under the assumption that
the sequence is strictly stationary but not necessarily independent. However, it was Kulik~(2008), who considered the
density deconvolution for LRD and short-range dependent (SRD) processes. However, Kulik~(2008) did not considered nonlinear wavelet estimators
but dealt instead with linear kernel estimators.\\

In nonparametric regression estimation, ARIMA-type models for the
errors  were analyzed in Cheng and Robinson~(1994), with error terms
of the form $\sigma (x_i, \xi_i)$. In Cs\"{o}rgo and
Mielniczuk~(2000), the error terms were modeled as infinite order
moving averages processes. Mielniczuk and Wu~(2004) investigated
another form of LRD, with the assumption that $x_i$ and $\xi_i$ are
not necessarily independent for the same $i$. ARIMA-type error
models were also considered in Kulik and Raimondo~(2009). In the
standard deconvolution model, and using a maxiset approach, Wishart
(2012) applied a fractional Brownian motion to model the presence of
LRD,  while Wang~(2012) used a minimax approach to study the problem
of recovering a function $f$ from a more general
noisy linear transformation where the noise is also a fractional Brownian motion.\\

The objective of this paper is to study the multichannel
deconvolution model from a minimax point of view, with the
relaxation that errors exhibit LRD. We do not limit our
consideration to a specific type of LRD: the only restriction is
that the errors should satisfy Assumption A1. In particular, we
derive minimax lower bounds for the $L^2$-risk in model \fr{sampl}
under  Assumption A1 when $f(\cdot)$ is assumed to belong to a Besov
ball and $g(\cdot,\cdot)$   has smoothness properties similar to
those in Pensky and Sapatinas~(2009, 2010), including both
regular-smooth and super-smooth convolutions. In addition, we
propose   an adaptive wavelet estimator for $f(\cdot)$ and  show
that such estimator is asymptotically optimal (or near-optimal
within a logarithmic factor) in the minimax sense, in a wide range
of Besov balls. We prove that the convergence rates of the resulting
estimators depend on the balance between the smoothness parameter
(of the response function $f(\cdot)$), the kernel parameters (of the
blurring function $g(\cdot,\cdot)$), and the long memory parameters
$d_l$, $l=1,2 \ldots, M$ (of the error sequence $\bxi^{(l)}$). Since
the parameters $d_l$ depend on the values of $l$, the convergence
rates have more complex expressions than the ones obtained in Kulik
and Raimondo~(2009) when studying nonparametric regression
estimation with ARIMA-type error models. The convergence rates we
derive are more similar in nature to those in Pensky and
Sapatinas~(2009, 2010). In particular, the convergence rates depend
on how the total number $n = N M$
of observations is distributed among the total number $M$ of channels. As we illustrate in two examples, convergence rates are not affected
by long range dependence in case of super-smooth convolutions, however, the situation changes in regular cases.  \\

The paper is organized as follows.  Section \ref{sec:long-range} discusses stationary sequences with LRD errors, justifies Assumption A1 and provides illustrative examples of stationary sequences satisfying this assumption.
Section \ref{sec:est_algorithm} describes the construction of the suggested wavelet estimator of $f(\cdot)$.
Section \ref{sec:lower_bounds} derives minimax lower bounds for the $L^2$-risk  for observations from model \fr{sampl}.
Section \ref{sec:upper_bounds}  proves that the suggested wavelet estimator is adaptive and asymptotically
optimal  (in the minimax sense) or near-optimal within a  logarithmic factor, in a wide range of Besov balls.  Section \ref{sec:examples} presents
examples of inverse problems in mathematical physics where one needs to recover initial
or boundary conditions on the basis of observations from a noisy solution of a partial differential equation to illustrate the application of the theory we developed.
Section \ref{sec:discussion} concludes with a brief discussion. Section \ref{sec:proofs} contains the
proofs of the theoretical results obtained in earlier sections.


 \section{Stationary Sequences with Long-Range Dependence}
\label{sec:long-range}
\setcounter{equation}{0}

In this section, for simplicity of exposition, we consider one sequence of errors  $\{ \xi_j:\, j=1,2,\ldots \}$.
Assume that $\{\xi_j:\, j=1,2,\ldots\}$ is a second-order stationary sequence
with covariance function $\gamma_{\xi}(k):=\gamma(k)$, $k=0,\pm1, \pm 2,\ldots$~. The spectral
density is defined as
$$
a_{\xi}(\lambda):=a(\lambda):=\frac{1}{2\pi}\sum_{k=-\infty}^{\infty}\gamma(k)
\exp(-i k\lambda)\;, \qquad \lambda\in [-\pi,\pi].
$$
On the other hand, the inverse transform which recovers $\gamma(k)$, $k=0,\pm1, \pm 2,\ldots$, from $a(\lambda)$, $\lambda\in [-\pi,\pi]$, is given by
$$
\gamma(k)=\int_{-\pi}^{\pi}e^{i k\lambda}a(\lambda)d\lambda, \quad k=0,\pm1, \pm 2,\ldots,
$$
under the assumption that the spectral density $a(\lambda)$, $\lambda\in [-\pi,\pi]$, is squared-integrable.\\

Let ${\bSigma}=[\gamma(j-k)]_{j,k=1}^N$ be the covariance
matrix of $(\xi_1,\ldots,\xi_N)$. Define ${\cal X}=\{\mathbf{x}\in
\CC^N:\mathbf{x}^*\mathbf{x}=1\}$, where $\mathbf{x}^*$ is
the complex-conjugate of $\mathbf{x}$. Since ${\bf \Sigma}$ is
Hermitian, one has
\begin{equation}\label{eq:spectral-radius}
\lambda_{\min} (\bSigma)=\inf_{\mathbf{x}\in{\cal X}}
\left(\mathbf{x}^*\bSigma \mathbf{x}\right) \qquad \text{and} \qquad \lambda_{\max}
(\bSigma)=\sup_{\mathbf{x}\in{\cal X}} \left(\mathbf{x}^*\bSigma
\mathbf{x}\right).
\end{equation}
With the definitions introduced above,
\begin{equation}\label{eq:spectral-repr}
\mathbf{x}^*\bSigma \mathbf{x}=\sum_{j,k=1}^N \mathbf{x}^*
\gamma(j-k) \mathbf{x}=\int_{-\pi}^{\pi}\left|\sum_{j=1}^N
x_je^{-ij\lambda}\right|^2 a(\lambda)d\lambda.
\end{equation}
Note that, by the Parseval identity, the function
$h(\lambda) = \left|\sum_{j=1}^N x_je^{-ij\lambda}\right|^2$, $\lambda \in [-\pi,\pi]$,
belongs to the set
$$
{\cal H}_N=\left\{h:\, h\; \mbox{\rm symmetric},\, |h|_{\infty}\le N,
\int_{-\pi}^{\pi}h(\lambda)d\lambda =2\pi \right\}.
$$
Let $d\in [0,1/2)$. Consider the following class of spectral densities
\be \label{eq:spec_den}
{\cal F}_d=\left\{a:\ a(\lambda)=|\lambda|^{-2d}a_{*}(\lambda), \;
 0<C_{\min}\le
|a_{*}(\lambda)| \leq C_{\max}<\infty,\,  \lambda\in [-\pi,\pi]\right\}.
\ee

Below we provide  two examples of second-order stationary sequences such that their spectral densities $a(\lambda)$, $\lambda \in [-\pi,\pi]$,
belong to the class ${\cal F}_d$ described in \fr{eq:spec_den}.
\\

{\bf Fractional ARIMA(0, $d$, 0).}
{\rm
Let $\{\xi_{j}:\, j=1,2,\ldots\}$ be the second-order stationary sequence
$$
\xi_{j}=\sum_{m=0}^{\infty}a_{m}\eta_{j-m},
$$
where $\eta_{j}$ are uncorrelated, zero-mean, random variables,
$\sigma_{\eta}^{2}:=\Var(\eta_{j})<\infty$, and
$$
a_{m}=(-1)^{m}\binom{-d}{m}=(-1)^{m}\frac{\Gamma(1-d)}{\Gamma(m+1)\Gamma(1-d-m)}
$$
with $d\in [0,1/2)$. Then, $a_{m}$, $m=0,1,\ldots$, are the coefficients in the power-series representation
$$
A(z):=(1-z)^{-d}:=\sum_{m=0}^{\infty}a_{m}z^{m}.
$$
Therefore, the spectral density $a(\lambda)$, $\lambda \in [-\pi,\pi]$, of $\{\xi_{j}:\, j=1,2,\ldots\}$, is given by
$$
a(\lambda) =\frac{\sigma_{\eta}^{2}}{2\pi}\left|
A(e^{-i\lambda})\right|^{2} =\frac{\sigma_{\eta}^{2}}{2\pi}
\left\vert 1-e^{-i\lambda}\right\vert ^{-2d}
=\frac{\sigma_{\eta}^{2}}{2\pi} \left\vert
2(1-\cos\lambda)\right\vert^{-d}
\sim\frac{\sigma_{\eta}^{2}}{2\pi}\left\vert \lambda\right\vert
^{-2d}\text{ }(\lambda\rightarrow 0).
$$
Hence, the sequence $\{\xi_{j}:\, j=1,2,\ldots\}$ has spectral density $a(\lambda)$, $\lambda \in [-\pi,\pi]$, that belongs to the class ${\cal F}_d$ described in \fr{eq:spec_den}. (The sequence $\{\xi_{j}:\, j=1,2,\ldots\}$ is
called the fractional ARIMA(0,$d$,0) time series.)
}
\\

{\bf Fractional Gaussian Noise.}
{\rm
Assume that $B_H(u)$, $u\in [0,\infty]$, is a fractional Brownian
motion with the Hurst parameter $H\in [1/2, 1)$. Define the second-order stationary sequence
$\xi_j=B_H(j)-B_H(j-1)$, $j=1,2,\dots$~. Its spectral density $a(\lambda)$, $\lambda \in [-\pi,\pi]$, is given by (see,
e.g., \cite{GPH}, p. 222)
$$
a(\lambda)=\sigma^2(2\pi)^{-2H-2}\Gamma(2H+1)\sin(\pi
H)4\sin^2(\lambda/2)\times
\sum_{k=-\infty}^{\infty}|k+(\lambda/2\pi)|^{-2H-1},
$$
and, hence,
$$
a(\lambda)=\frac{2\sigma^2}{\pi}\Gamma(2H+1)\sin(\pi H)
\lambda^{1-2H}\; \qquad (\lambda \to 0).
$$
Hence, the sequence $\{\xi_{j}:\, j=1,2,\ldots\}$ has spectral density $a(\lambda)$, $\lambda \in [-\pi,\pi]$, that belongs to class ${\cal F}_d$ with
$d=H-1/2$. (The sequence $\{\xi_{j}:\, j=1,2,\ldots\}$ is called the
fractional Gaussian noise.)
}
\\

It follows from \fr{eq:spec_den} that, for $a\in {\cal F}_d$, one has $a(\lambda) \sim
|\lambda|^{-2d}$ ($\lambda\to 0$).
It also  turns out that the condition $a\in {\cal F}_d$, $d\in [0,1/2)$,
implies that all eigenvalues of the covariance matrix $\bSigma$ are of asymptotic order $N^{2d}$ ($N \to \infty$). In particular, the following lemma is true.

\begin{lemma}\label{lemma:assump1}
Assume that $\{\xi_{j}:\, j=1,2,\ldots\}$ is a second-order stationary sequence
with spectral density $a\in {\cal F}_d$, $d\in [0,1/2)$. Then,
for some constants $K_{1 d}$ and $K_{2 d}$ $(0 < K_{1 d} \leq K_{2 d} < \infty)$, that depend on $d$ only,
$$
K_{1 d} N^{2d} \leq \lambda_{\min} (\bSigma) \leq \lambda_{\max} (\bSigma) \leq K_{2 d} N^{2d}.
$$
\end{lemma}

\begin{remark}{\rm
If $d=0$, then ${\cal F}_d$ is the class of spectral
densities $a(\lambda)$ that are bounded away from 0 and $\infty$ for all $\lambda \in [-\pi,\pi]$.
In particular, the corresponding second-order stationary sequences $\{\xi_{j}:\, j=1,2,\ldots\}$ are weakly dependent.
Then, the statement of Lemma \ref{lemma:assump1} reduces to a result
in Grenander and Szeg\"{o} \cite{GSz}, Section 5.2.}
\end{remark}

It follows immediately from Lemma \ref{lemma:assump1} that if, for each $l=1,2,\ldots,M$, $\bxi^{(l)}$ is a second-order stationary Gaussian sequence
with spectral density $a_l \in {\cal F}_{d_l}$, $d_l \in [0,1/2)$, that $\bxi^{(l)}$ are independent for different $l$'s, and that $d_l$'s are uniformly bounded, then
Assumption A1 holds.

\begin{corollary} \label{cor1}
For each\, $l=1,2,\ldots,M$, let $\bxi^{(l)}$ be a second-order stationary Gaussian sequence
with spectral density $a_l \in {\cal F}_{d_l}$, $d_l \in [0,1/2)$. We assume that $\bxi^{(l)}$ are independent for different $l$'s.
Let $d_l$, $l=1,2,\ldots,M$, be uniformly bounded, i.e., there exists $\dus$ $(0 \leq \dus <1/2)$ such that, for each $l=1,2,\ldots,M$,
\beqn \label{eq:d_cond}
 0 \leq d_l \leq \dus < 1/2.
\eeqn
Then, Assumption A1 holds.
\end{corollary}


 \section{The Estimation Algorithm}
\label{sec:est_algorithm}
\setcounter{equation}{0}

In what follows, $\langle \cdot,\cdot \rangle$ denotes the inner
product in $\RR^N$.  We also denote the complex-conjugate of $a \in \CC$ by $\bar{a}$,
 the discrete Fourier basis on  the interval  $T$  by $e_m(t_i) = e^{-i 2 \pi m t_i}$, $t_i=i/N$, $i=1,2, \ldots, N$,
$m=0,\pm1, \pm 2,\ldots$, and the complex-conjugate of the matrix $\bA$ by $\bA^*$. \\

Recall the multichannel deconvolution model \fr{sampl}. Denote
\begin{eqnarray*}
h(u_l, t_i)= \int_{T}  g(u_l, t_i -x) f( x) dx, \quad l=1,2,\ldots,M, \quad i=1,2,\ldots,N.
\end{eqnarray*}
Then, equation \fr{sampl} can be rewritten as
\be \label{samplh}
y(u_l, t_i) = h(u_l, t_i)  +  \xi_{li}, \quad l=1,2,\ldots,M, \quad i=1,2,\ldots,N.
\ee
For each $l=1,2,\ldots,M$, let   $h_m(u_l) = \langle e_m, h(u_l, \cdot) \rangle$,
$y_m(u_l) = \langle e_m, y(u_l, \cdot) \rangle$, $z_{lm}= \langle e_m, \bxi^{(l)} \rangle$,
$g_m(u_l)= \langle e_m, g(u_l, \cdot) \rangle$ and $f_m= \langle e_m, f \rangle$
be the discrete Fourier coefficients of the $\RR^N$ vectors $h(u_l, t_i)$, $y(u_l, t_i)$, $\xi_{li}$, $g(u_l, t_i)$ and $f(t_i)$,  $i=1,2,\ldots,N$, respectively.
Then, applying the discrete Fourier transform to
\fr{samplh}, one obtains, for any $u_l\in U$,  $l=1,2,\ldots,M$,
\be \label{hgf2}
y_m(u_l)= g_m(u_l)f_m + N^{-1/2} z_{lm}
\ee
and
\be \label{hgf}
h_m(u_l)= g_m(u_l) f_m.
\ee

Multiplying both sides of \fr{hgf2} by $ N^{-2 d_l}  \overline{g_m(u_l)}$, and adding them together, we obtain the following estimator of $f_m$
\be \label{fmh}
\widehat{f}_m= \left( \sum^M_{l=1}  N^{-2 d_l}\,  \overline{g_m(u_l)} y_m(u_l) \right)
/\left( \sum^{M}_{l=1}N^{-2 d_l}  |g_m(u_l)|^2 \right).
\ee
Let  $\ph^*(\cdot)$ and $\psi^*(\cdot)$ be the Meyer scaling and
mother wavelet functions, respectively, defined on the real line (see, e.g., Meyer (1992)),
and  obtain a periodized version of Meyer
wavelet basis as in Johnstone\ {\it et al.} (2004), i.e., for $j\geq 0$ and $k=0,1,\ldots,2^j-1$,
$$
\ph_{jk}(x) = \sum_{i \in \ints} 2^{j/2} \ph^*(2^j (x +i) -  k),
\quad \psi_{jk}(x) = \sum_{i \in \ints} 2^{j/2} \psi^*(2^j (x +i) -
k), \quad x \in T.
$$
Following Pensky and Sapatinas (2009, 2010), using the periodized Meyer wavelet basis described above,
for some $j_0 \geq 0$, expand $f(\cdot) \in L^2(T)$  as
\be
f(t) = \sum_{k=0}^{2^{j_0}-1} a_{j_0k} \ph_{j_0k}
(t) + \sum_{j=j_0}^\infty \sum_{k=0}^{2^j -1} b_{jk} \psijk (t), \quad t \in T.
\label{funf}
\ee
Furthermore, by Plancherel's formula, the scaling
coefficients, $\ajk=\langle f, \ph_{j_0k} \rangle$, and the wavelet
coefficients, $\bjk=\langle f,\psi_{jk}\rangle$, of $f(\cdot)$ can
be represented as
\be  \label{alkandblk}
 \ajk = \sum_{m \in C_{j_0}} \fm \overline{\phmjk},
\ \ \ \bjk = \sum_{m \in C_j} \fm \overline{\psimjk},
\ee
where  $C_{j_0} = \lfi m: \ph_{mj_0 k} \neq 0 \rfi$ and, for any $j \geq j_0$,
$$
C_j = \lfi m: \psimjk \neq 0 \rfi  \subseteq 2\pi/3 [-2^{j+2}, -2^j] \cup [2^j, 2^{j+2}].
$$
(Note that the cardinality  $|C_j|$ of the set $C_j$ is $|C_j|
= 4 \pi 2^j $, see, e.g., Johnstone\ {\it et al.} (2004).)
Estimates of $\ajk$ and $\bjk$ are readily obtained by substituting $\fm$ in \fr{alkandblk}
with   \fr{fmh}, i.e.,
\be \label{coefest}
\hajk = \sum_{m \in C_{j_0}} \hfm \overline{\phmjk},\ \ \ \hbjk = \sum_{m \in C_j} \hfm
\overline{\psimjk}.
\ee

We now construct a (block thresholding) wavelet estimator of
$f(\cdot)$, suggested by Pensky \& Sapatinas (2009, 2010). For this
purpose, we divide the wavelet coefficients at each resolution level
into blocks of length $\ln n$. Let $\Aj$ and $\Ujr$ be the following
sets of indices
$$
\Aj = \lfi r \mid   r=1,2,\ldots, 2^j/\ln n \rfi,
$$
$$ \Ujr = \lfi
k \mid k = 0,1, \ldots, 2^j-1;\ (r-1) \ln n \leq k \leq r \ln n -1
\rfi.
$$
Denote
\be \label{bjr}
\Bjr = \sumku \bjk^2,\ \ \ \hBjr = \sumku \hbjk^2.
\ee
Finally, for any $j_0 \geq 0$, the (block thresholding) wavelet estimator $\hat{f}_n(\cdot)$ of $f(\cdot)$ is constructed as
\be  \label{fest}
\hat{f}_n(t) = \sum_{k=0}^{2^{j_0} -1} \hajk
\ph_{j_0k} (t) + \sum_{j=j_0}^{J-1} \sumr \sumku  \hbjk \II(|\hBjr|
\geq \lam_{j})\, \psijk (t), \quad t \in T,
\ee
where $\II(A)$ is the indicator function of the set $A$, and the resolution levels $j_0$
and $J$ and the thresholds $\lam_{j}$ will be defined in Section~\ref{sec:upper_bounds}.

\medskip

In what follows, the symbol $C$ is used for a generic positive
constant, independent of $n$, while the symbol $K$ is used for a generic positive
constant, independent of $m$, $n$, $M$ and $u_1,u_2,\ldots,u_M$. Either of $C$ or $K$ may take different values at different places.


 \section{ Minimax Lower Bounds for the $L^2$-Risk}
\label{sec:lower_bounds}
 \setcounter{equation}{0}

 Denote
 \be \label{sdef}
s'= s + 1/2 -1/p, \quad s^*= s + 1/2 -1/p', \quad p'=\min\{ p, 2\}.
\ee Assume that the unknown response function $f(\cdot)$ belongs to
a Besov ball $B^{s}_{p, q}(A)$ of radius $A>0$, so that   the
wavelet coefficients $a_{j_0k}$ and $b_{jk}$ defined in
\fr{alkandblk} satisfy the following relation \be \label{assum1}
B^{s}_{p, q}(A)= \left\{ f  \in L^2(U):\ \left( \sum^{2^{j_0}-1}_{k
= 0} \left| a_{j_{0}k}\right|^p \right)^{\frac{1}{p}}+  \left(
\sum^{\infty}_{j =j_0} 2^{j s'q} \left( \sum^{2^j-1}_{k=0} \left|
b_{jk}\right|^p \right)^{\frac{q}{p}} \right)^{1/q} \leq A \right\}.
\ee Below, we construct  minimax lower bounds for the (quadratic)
$L^2$-risk. For this purpose, we define the minimax $L^2$-risk over
the set $V \subseteq L^2(T)$ as
$$
R_{n} (V) = \inf_{\tilde{f}} \, \sup_{f \in V} \EE \| \tilde{f}  - f \|^2,
$$
where $\| g \|$ is the $L^2$-norm of a function $g(\cdot)$ and the
infimum is taken over all possible estimators $\tilde{f}(\cdot)$
(measurable functions taking their values in a set containing $V$) of $f(\cdot)$,
based on observations from model \fr{sampl}).\\

For $M = M_n$ and $N = n/M_n$, denote
\be  \label{tau_def}
\tau_{\kappa} (m,n)=  M^{-1} \sum^M_{l=1}  N^{-2 \kappa d_l}\, |g_m(u_l)|^{2\kappa}, \quad \kappa=1 \; \text{or} \; 2 \; \text{or} \; 4,
\ee
and
\be \label{Delt}
\Delta_\kappa (j,n) =
 |C_j|^{-1}\, \sum_{m \in C_j} \tau_{\kappa} (m, n)\ [\tau_1 (m, n)]^{-2 \kappa},\ \ \kappa =1 \; \text{or} \; 2.
\ee

The expression $\tau_1 (m,n)$ appears in both the lower and the upper bounds for the $L^2$-risk. Hence, we
impose the following assumption:
\\

{\bf Assumption A2}:\ \  For some  constants $\nu_1, \nu_2, \lambda_1, \lambda_2  \in \RR$,
$\alpha_1, \alpha_2 \geq 0$ ($\lambda_1, \lambda_2 >0$ if $\alpha_1 = \alpha_2 = 0$, $\nu_1 = \nu_2 =0$)\ and $K_3, K_4, \beta > 0$,  independent of $m$ and $n$, and for some sequence $\eps_n>0$,
independent of $m$,  one has
\be \label{eq:assump2}
K_3 \eps_n\ |m|^{-2\nu_1}  (\ln|m|)^{-\lambda_1} e^{ -\alpha_1 |m|^\beta} \leq  \tau_1 (m, n)
\leq K_4 \eps_n\ |m|^{-2\nu_2}  (\ln|m|)^{-\lambda_2} e^{ -\alpha_2 |m|^\beta},
\ee
where   either $\alpha_1 \alpha_2 \neq 0$ or $\alpha_1 = \alpha_2 = 0$ and $\nu_1 = \nu_2 = \nu >0$.
The sequence $\eps_n$  in \fr{eq:assump2} is such that
\be \label{ns_prop1}
\ns = n \eps_n \rightarrow \infty \quad (n \rightarrow \infty).
\ee

Under Assumptions A1 and A2, the following statement is true.

\begin{theorem}    \label{th:lower}
Let Assumptions A1 and A2  hold. Let $\{\phi_{j_0,k}(\cdot),\psi_{j,k}(\cdot)\}$ be the periodic Meyer wavelet basis discussed in Section
\ref{sec:est_algorithm}.  Let $s >\max(0,1/p-1/2)$, $1 \leq p \leq \infty$, $1 \leq q \leq \infty$ and
$A>0$.  Then, as $n \rightarrow \infty$,
\be
R_n  (B_{p,q}^s (A)) \geq \lfi
\begin{array}{ll}
C (\ns)^{-\frac{2s}{2s+2\nu+1}}\ (\ln \ns)^{\frac{2s \lambda_2}{2s+2\nu+1}},
& {\rm if}\;\;\; \alpha_1 = \alpha_2 =0,\  \nu(2-p) < p{s^*},
\\
C \Big ( \frac{\ln \ns}{\ns} \Big )^{\frac{2{s^*}}{2s^*+2\nu}}
(\ln \ns)^{\frac{2{s^* \lambda_2}}{2s^*+2\nu}}, & {\rm
if}\;\;\; \alpha_1 = \alpha_2 =0,\   \nu(2-p) \geq p{s^*},
\label{low1} \\
C (\ln \ns)^{-\frac{2{s^*}}{\beta}}, &  \rm{if}\,\;\; \alpha_1 \alpha_2 \neq 0.
\end{array} \right.
\ee
\end{theorem}


\section {Minimax Upper Bounds for the $L^2$-Risk}
 \label{sec:upper_bounds}
\setcounter{equation}{0}

Let $\hat{f}_n(\cdot)$ be the (block thresholding) wavelet estimator defined by \fr{fest}.
Choose now $\jo$ and $J$ such that
\beqn
2^\jo = \ln \ns, \ \ 2^J = (\ns)^{\frac{1}{2\nu+1}},& \mbox{if}\;\;\; \alpha_1=\alpha_2=0,  \label{jpower} \\
2^\jo = \frac{3}{8 \pi} \lkr \frac{\ln \ns}{2\alpha}
\rkr^{\frac{1}{\beta}}, \ \ 2^J=2^{\jo}, & \mbox{if}\;\;\; \alpha_1 \alpha >0.
\label{jexp}
\eeqn
(Since $\jo >J-1$ when $\alpha_1 \alpha >0$, the estimator
\fr{fest} only consists of the first (linear) part and, hence,
$\lam_j$ does not need to be selected in this case.) Set, for some
constant $\mu>0$, large enough,
\be
\lam_j = \mu^2\; (\ns)^{-1}  \ln (\ns)\, 2^{2 \nu j}  j^{\lam_1}, \ \ \
\mbox{if}\ \ \alpha_1=\alpha_2=0.
\label{lamj}
\ee
Note that the choices of $\jo$,  $J$ and $\lam_j$ are independent of the parameters, $s$,
$p$, $q$ and $A$ of the Besov ball $B_{p,q}^s (A)$; hence, the estimator \fr{fest} is
adaptive with respect to these parameters.\\

Denote $(x)_{+} = \max(0,x)$,
\be \varrho = \lfi
\begin{array}{ll}
\frac{(2\nu +1)(2-p)_+}{p(2s+2\nu+1)}, & \mbox{if}\;\;\; \nu(2-p) <
p{s^*}, \\
\frac{(q-p)_+}{q}, & \mbox{if}\;\;\; \nu(2-p) = p{s^*}, \label{rovalue}\\
0, & \mbox{if}\;\;\; \nu(2-p) > p{s^*}.  
\end{array} \right.
\ee
Assume that,
in the case of $\alpha_1 = \alpha_2 =0$, the sequence $\eps_n$  is  such that
\be \label{ns_prop2}
-h_1 \ln n \leq \ln(\eps_n) \leq  h_2  \ln n
\ee
for some constants   $h_1, h_2 \in (0,1)$. Observe that condition \fr{ns_prop2} implies
\fr{ns_prop1} and that $\ln \ns \asymp \ln n$ ($n \rightarrow \infty$).
(Here, and in what follows, $u(n) \asymp v(n)$ means that there exist constants $C_1, C_2$ $(0 < C_1 \leq C_2 < \infty)$, independent of $n$, such that $0<C_1 v(n) \leq u(n) \leq C_2 v(n)<\infty$ for $n$ large enough.)
\\


The proof of the minimax upper bounds for the $L^2$-risk is based on
the following two lemmas.

\begin{lemma} \label{l:coef}
Let Assumptions A1 and A2 hold. Let the estimators
$\hajk$ and $\hbjk$ of the scaling and wavelet coefficients $\ajk$
and $\bjk$, respectively, be given by \fr{alkandblk} with
$\hfm$ defined by \fr{fmh}. Then,  for
all $j \geq j_0$,
\be  \label{ha}
\EE |\hajk - \ajk|^2   \leq  C n^{-1} \Delta_1 (\jo,n)\quad \text{and} \quad
\EE |\hbjk - \bjk|^{2}   \leq  C n^{-1} \Delta_1 (j,n).
\ee
If $\alpha_1=\alpha_2=0$ and \fr{ns_prop2} holds, then, for any $j \geq j_0$,
\begin{align}
\EE |\hbjk - \bjk|^{4}  & \leq  C  n^3\, (\ln n)^{3 \lam_1} \, (\ns)^{-\frac{3}{2 \nu +1}}.
\label{hbb}
\end{align}
 \end{lemma}

\begin{lemma} \label{l:deviation}
Let Assumptions A1, A2  and \fr{ns_prop2} hold. Let the estimators $\hbjk$ of the wavelet coefficients $\bjk$ be
given by  \fr{alkandblk} with $\hfm$ defined by
\fr{fmh}. Let
\be \label{mu_cond}
\mu \geq \sqrt{\frac{2}{1-h_1}}\ \lkv \sqrt{c_1} + \frac{\sqrt{8 \pi \kappa}}{\sqrt{K_3}} (\ln 2)^{\lambda_1/2} \lkr\frac{2\pi}{3} \rkr^\nu \rkv,
\ee
where $c_1$, $K_3$ and $h_1$ are defined in \fr{delta1upper}, \fr{eq:assump2} and \fr{ns_prop2}, respectively. Then,
 for all $j \geq j_0$ and any $\kappa >0$,
\be \label{probbound}
\PP \lkr  \sumku |\hbjk - \bjk|^2 \geq  (4 \ns)^{-1}\ \mu^2\,  2^{2\nu j}\,  j^\lam\, \ln \ns
  \rkr \leq  n^{-\kappa}.
\ee
\end{lemma}

Under Assumptions A1 and A2, and using Lemmas \ref{l:coef} and \ref{l:deviation}, the following statement is true.

\begin{theorem}  \label{th:upper}
Let Assumptions  A1 and A2 hold. Let $\hat{f}_n(\cdot)$ be the
wavelet estimator defined by \fr{fest}, with $\jo$ and $J$ given by
\fr{jpower} (if $\alpha_1 = \alpha_2 = 0$) or \fr{jexp} (if $\alpha_1 \alpha_2 >0$) and $\mu$ satisfying \fr{mu_cond} with $\kappa =5$.
Let $s > 1/p'$, $1 \leq p \leq \infty$, $1 \leq q \leq \infty$ and $A>0$.
Then, under \fr{ns_prop1} if $\alpha_1 \alpha_2 >0$ or  \fr{ns_prop2}
if $\alpha_1 = \alpha_2 =0$, as $n \rightarrow \infty$,
\be
\sup_{f \in B_{p,q}^s (A)} \EE \|\hat{f}_n -f\|^2 \leq \lfi
\begin{array}{ll}
C  (\ns)^{-\frac{2s}{2s+2\nu+1}}\ \lkr  \ln n\rkr^{\varrho + \frac{2s\lambda_1}{2s+2\nu +1}}, & {\rm if}\;\;\;
\alpha_1 = \alpha_2 =0,\ \nu(2-p) < p{s^*},
\\ 
C \Big ( \frac{\ln n}{\ns} \Big )^{\frac{2{s^*}}{2s^*+2\nu}} \ \lkr
\ln n \rkr^{\varrho+\frac{2{s^* \lambda_1 }}{2s^*+2\nu }}, & {\rm if}\;\;\; \alpha_1 = \alpha_2 =0,\ \nu(2-p) \geq p{s^*},
\label{up} \\
C (\ln  \ns )^{-\frac{2{s^*}}{\beta}}, & {\rm if}\;\;\; \alpha_1 \alpha_2 >0.
\end{array} \right.
\ee
\end{theorem}

\begin{remark}
{\rm Theorems \ref{th:lower} and \ref{th:upper} implies that, for the $L^2$-risk, the wavelet estimator $\hat{f}_n(\cdot)$ 
defined by \fr{fest} is asymptotical optimal (in the minimax sense), or near optimal within a logarithmic factor, over a wide 
range of Besov balls $B_{p,q}^s(A)$ of radius $A>0$ with $s>\max(1/p,1/2)$, $1 \leq p \leq \infty$ and $1 \leq q \leq \infty$. 
The convergence rates depend on the balance between the smoothness parameter $s$
(of the response function $f(\cdot)$), the kernel parameters $\nu, \beta, \lambda_1$ and $\lambda_2$ (of the blurring function 
$g(\cdot,\cdot)$), the long memory parameters $d_l$, $l=1,2 \ldots, M$ (of the error sequence $\bxi^{(l)}$), 
and how the total number of observations $n$ is distributed among the total number of channels $M$.  In particular, $M$ and $d_l$, $l=1,2, \ldots, M$, jointly determine the value of $\eps_n$ which, in turn,
defines the ``essential'' convergence rate $n^* = n \eps _n$ which may differ considerably from $n$. 
For example, if $M= M_n = n^\te$, $0 \leq \te <1$ and $|g_m(u_l)|^{2} \asymp |m|^{- 2 \nu}$ for every $l =1,2 \ldots, M$,
then  
\be  \label{eq:epsn_expr}
\eps_n =  M^{-1} \sum_{l=1}^M N^{-2 d_l},
\ee
and, therefore, $n^{1 - 2 \dus (1 - \te)} \leq \ns \leq n$, where $d^* = \max_{1\leq l\leq M} d_l$, so that, 
 $\ns$ can take any value between    $n^{1 - 2 \dus (1 - \te)}$ and $n$. This is further illustrated in Section \ref{sec:examples} below.
%
%
}
\end{remark}


\section {Illustrative Examples}
 \label{sec:examples}
\setcounter{equation}{0}

In this section, we consider some illustrative examples of application of the theory developed in the previous sections. They are particular examples
of inverse problems in mathematical physics where one needs to recover initial
or boundary conditions on the basis of observations from a noisy solution of a partial differential equation.\\

We assume that condition \fr{eq:d_cond} holds true and that there exist   $\te_1$ and $\te_2$, such that $M=M_n$ satisfies
\beqn
 n^{\te_1} \leq M \leq n^{\te_2}, \quad 0 \leq \te_1 \leq \te_2 <1. \label{eq:M_cond}
\eeqn
(Note that, under (\ref{eq:M_cond}), $n^{1 - \te_2} \leq N \leq n^{1 - \te_1}$.)


\begin{example}
\label{ex1}
{\rm Consider the case when $g_m(\cdot)$, $m=0,\pm1, \pm 2,\ldots$, is of the form
\be \label{gm_ex1}
g_m(u) = C_g \exp \lkr - K |m|^\beta q(u) \rkr, \quad u \in U,
\ee
where $q(\cdot)$ in \fr{gm_ex1} is such that, for some $q_1$ and $q_2$,
\be \label{eq:q_cond}
0  < q_1 \leq q(u) \leq q_2 < \infty, \quad  u \in U.
\ee

This set up takes place in the estimation of the initial condition in the heat
conductivity equation or the estimation of the boundary condition for the
Dirichlet problem of the Laplacian on the unit circle (see Examples 1 and 2 of Pensky and Sapatinas (2009, 2010)). In the former case,
 $g_m(u) = \exp(- 4 \pi^2 m^2 u)$, $u \in U$, so that $K = 4 \pi^2$,
$\beta = 2$, $q(u) = u$, $q_1 =a$ and $q_2 =b$. In the latter case,
$g_m(u) = C u^{|m|} = C \exp (- |m| \ln(1/u))$, $0< r_1 \leq u \leq r_2 < 1$,  so that
$K=1$, $\beta=1$, $q(u) = \ln(1/u)$, $q_1 = \ln(1/r_2)$ and  $q_2 = \ln(1/r_1)$.\\

It is easy to see that, under conditions \fr{gm_ex1} and \fr{eq:q_cond},
 for $\tau_1 (m,n)$ given in \fr{tau_def},
$$
\tau_1(m,n) \leq  C_g \, \eps_n\, \exp \lkr - 2 K q_1 |m|^\beta  \rkr \quad \text{and} \quad
\tau_1(m,n) \geq  C_g \, \eps_n\, \exp \lkr - 2 K q_2 |m|^\beta  \rkr,
$$
where $\eps_n$ is of the form \fr{eq:epsn_expr}.
Assumptions \fr{eq:d_cond} and \fr{eq:M_cond} lead to the following bounds for $\ns$:
$$
n^{1 - 2 \dus (1 - \te_1)} \leq \ns \leq n,
$$
so that $\ln n \asymp \ln \ns$. Therefore, according to Theorems \ref{th:lower} and \ref{th:upper},
\be  \label{rn_ex1}
R_n  (B_{p,q}^s (A)) \asymp  (\ln n)^{-\frac{2{s^*}}{\beta}}.
\ee

Note that, in this   case,  the value of $\dus$ has absolutely no
bearing on the convergence rates of the linear wavelet estimators:
the convergence rates  are determined entirely by the properties of
the smoothness parameter $s^{*}$ (of the response function
$f(\cdot)$) and the kernel parameter $\beta$ (of the blurring
function $g(\cdot,\cdot)$).\\

In other words, in case of super-smooth convolutions, LRD does not influence the convergence rates of the suugested wavelet
estimator. A similar effect is observed in the case of kernel smoothing,
see Section 2.2 in Kulik (2008).
 }

\end{example}


\begin{example}
\label{ex2}
{\rm

Suppose that the blurring function $g(\cdot,\cdot)$ is of a box-car like kernel,
i.e.,
\begin{equation}
\label{eq:box-car1}
g(u,t) = 0.5\, q(u)\, \II(|t| < u),\ \ \
u \in U, \;\; t \in T,
\end{equation}
where $q(\cdot)$ is some positive function which satisfies conditions \fr{eq:q_cond}. In this case,
the functional Fourier coefficients $g_m(\cdot)$ are of the form
\be
\label{boxcartype}
g_0(u)=1 \quad \text{and} \quad g_m(u)= (2 \pi m)^{-1}\, \gamma(u) \sin(2\pi m u),
\quad m \in \ints \setminus \{0\}, \quad \quad u \in U.
\ee

It is easy to see that estimation of the initial speed of a wave on a finite interval
(see Example 4 of Pensky and Sapatinas (2009) or Example 3 of Pensky and Sapatinas (2010)) leads to $g_m (\cdot)$ of the
form \fr{boxcartype} with $q(u) =1$.  Assume, without loss of generality, that $u \in [0,1]$,
so that $a=0$, $b=1$, and consider (equispaced channels) $u_l = l/M$, $l=1,2,\ldots,M$, such that
\be \label{eq:linear_d}
d_l = a_1 u_l + a_2, \quad 0 \leq a_2  \leq \dus < 1/2, \ 0 \leq a_1 + a_2 \leq \dus < 1/2,
\ee
i.e., condition \fr{eq:d_cond}  holds.  Note that if  $a_1 =0$, then 
$$
\tau_1 (m,n) \asymp  M^{-1} N^{-2a_2} (4 \pi^2 m^2)^{-1} \,  \sum_{l=1}^M  \sin^2(2 \pi m l /M),
$$
which is similar to the expression for $\tau_1 (m,n)$  studied in Section 6 of  Pensky and Sapatinas (2010).
Following their calculations, one obtains that, if $j_0$ in \fr{fest}
is such that $2^{j_0} > (\ln n)^\delta$ for some $\delta >0$ and
$M \geq (32 \pi/3)  n^{1/3}$, then, for $n$ and $|m|$ large enough,
$$ \tau_1 (m, n) \asymp  N^{-2a_2} m^{-2}.$$

Assume now, without loss of generality, that  $a_1 \geq 0$. (Note that the case of $a_1 \leq 0$
can be handled similarly by changing $u$ to $1-u$.) Below, we shall show that, in this case, a 
similar result can be obtained under less stringent conditions on $M = M_n$.
Indeed, the following statement is true. \\

\begin{lemma}  \label{lem:boxcar}
Let $g(\cdot,\cdot)$ be of the form
\fr{eq:box-car1}, where $q(\cdot)$ is some positive function which satisfies \fr{eq:q_cond},
and let $d_l$, $l=1,2,\ldots,M$, be given by \fr{eq:linear_d} with  $a_1 \geq 0$. Assume (without loss of generality) that  
$U=[0,1]$, and consider $u_l = l/M$, $l=1,2,\ldots, M$. Let  $M = M_n$ satisfy  \fr{eq:M_cond} with $\te_1>0$ if $a_1 >0$ and 
$M \geq (32 \pi/3)  n^{1/3}$ if $a_1 =0$. 
If  $m \in A_j$, where $|A_j|= C_m 2^j$, for some absolute constant $C_m > 0$,
with $j \geq j_0 >0$, where $j_0$ is such that $2^{j_0} \geq C_0  \ln n$  for some $C_0 > 0$,
then, for $n$ and $|m|$ large enough,
\begin{equation}
\label{eq:boxcar} \tau_1 (m,n) \asymp   N^{-2a_2} m^{-2} (\log n)^{-1}.
\end{equation}
\end{lemma}

It follows immediately from Lemma \ref{lem:boxcar} that, if
$$
M=M_n = n^{\te}, \quad 0< \te  <1,
$$
then  Assumption A2 holds with $\alpha_1 = \alpha_2 = 0$,
$\nu_1 = \nu_2 = \nu = 2$, $\eps_n = n^{-2 a_2 (1 - \te)}\, (\ln n)^{-1} $ and
$\lambda_1 = \lambda_2 = 0$. Note that $\eps_n$ satisfies
conditions \fr{ns_prop1} and \fr{ns_prop2}, so that $\ln n \asymp \ln \ns$.
Therefore, according to
Theorems \ref{th:lower} and \ref{th:upper},
\be  \label{lowbox}
R_n  (B_{p,q}^s (A)) \geq \lfi
\begin{array}{ll}
C (\ns)^{-\frac{2s}{2s+5}},
& {\rm if}\;\;\;    4-2p < p{s^*},
\\
C \lkr \frac{\ln \ns}{\ns} \rkr^{\frac{s^*}{s^*+2}},
& {\rm if}\;\;\;   4-2p \geq p{s^*},
\end{array} \right.
\ee
and
\be  \label{upperbox}
\sup_{f \in B_{p,q}^s (A)} \EE \|\hat{f}_n -f\|^2 \leq \lfi
\begin{array}{ll}
C  (\ns)^{-\frac{2s}{2s+5}}\ \lkr  \ln n \rkr^{\varrho},
& {\rm if}\;\;\;    4-2p < p{s^*},\\
C \Big ( \frac{\ln n}{\ns} \Big )^{\frac{s^*}{s^*+2}} \
\lkr \ln n \rkr^{\varrho},
& {\rm if}\;\;\;   4-2p \geq p{s^*},
\end{array} \right.
\ee
where $$\ns = n^{1 - 2 a_2 (1 - \te)}\, (\ln n)^{-1}$$  and
$$\varrho = \lfi
\begin{array}{ll}
\frac{(5(2-p)_+}{p(2s+5)}, & \mbox{if}\;\;\; 4-2p < p{s^*}, \\
\frac{(q-p)_+}{q}, & \mbox{if}\;\;\; 4-2p = p{s^*},  \\
0, & \mbox{if}\;\;\; 4-2p > p{s^*}.  
\end{array} \right.
$$
\\

Note that LRD affects the convergence rates in this case via the
parameter $a_2$ that appears in the definition (\ref{eq:linear_d}).

}

\end{example}


\section{Discussion. }
\label{sec:discussion} \setcounter{equation}{0} Deconvolution is the
common problem in many areas of signal and image processing which
include, for instance, LIDAR (Light Detection and Ranging) remote
sensing and reconstruction of blurred images. LIDAR is a laser
device which emits pulses, reflections of which are gathered by a
telescope aligned with the laser (see, e.g., Park, Dho \& Kong
(1997) and Harsdorf \& Reuter (2000)). The return signal is used to
determine distance and the position of the reflecting material.
However, if the system response function of the LIDAR is longer than
the time resolution interval, then the measured LIDAR signal is
blurred and the effective accuracy of the LIDAR decreases. If $M$
($M\geq2$) LIDAR devices are used to recover a signal, then we talk
about a multichannel deconvolution
problem. This leads to the discrete model (\ref{sampl}) considered in this work.\\

The multichannel deconvolution model (\ref{sampl}) can also be thought of as the discrete version of a model referred to as the functional deconvolution model by Pensky and Sapatinas (2009, 2010). The functional deconvolution model has a multitude of applications. In particular, it can be used in a number of inverse problems in mathematical physics where one needs to recover initial
or boundary conditions on the basis of observations from a noisy solution of a partial differential equation.  Lattes \& Lions
(1967) initiated research in the problem of recovering the initial
condition for parabolic equations based on observations in a
fixed-time strip. This problem and the problem of recovering the
boundary condition for elliptic equations based on observations in
an internal domain were studied in Golubev \& Khasminskii (1999);
the latter problem was also discussed in Golubev (2004). Some of these specific models were considered in Section  \ref{sec:examples}.\\

The multichannel deconvolution model (\ref{sampl}) and its continuous version, the functional deconvolution model, were studied by Pensky and Sapatinas~(2009, 2010), under the assumption  that errors are independent and identically distributed Gaussian random variables. The objective of this work was to study the multichannel deconvolution model (\ref{sampl}) from a minimax point of view, with the relaxation that errors exhibit LRD. We were not limited our consideration to a specific type of LRD: the only restriction made was that the errors should satisfy a general assumption in terms of the smallest and larger eigenvalues of their covariance matrices. In particular, minimax lower bounds for the $L^2$-risk in model \fr{sampl} under  such assumption were derived when $f(\cdot)$
is assumed to belong to a Besov ball and $g(\cdot,\cdot)$  has smoothness properties similar to those in Pensky and Sapatinas~(2009, 2010), including both regular-smooth and super-smooth convolutions.
In addition, an adaptive wavelet estimator of $f(\cdot)$ was constructed and shown that such estimator is asymptotically optimal
(in the minimax sense), or near-optimal within a logarithmic factor, in a wide range of Besov balls.  The convergence rates of the resulting estimators depend on the balance between the smoothness parameter
(of the response function $f(\cdot)$), the kernel parameters (of the blurring function $g(\cdot,\cdot)$), and the long memory parameters $d_l$, $l=1,2 \ldots, M$ (of the error sequence $\bxi^{(l)}$), and how the total number of observations is distributed among the total number of channels.
Note that SRD is implicitly included in our results by selecting $d_l=0$, $l=1,2,\ldots,M$. In this case, the convergence rates we obtained coincide with the convergence rates obtained under the assumption of independent and identically distributed Gaussian errors by Pensky and Sapatinas (2009, 2010).
\\

Under the assumption that the errors are independent and identically distributed Gaussian random variables, for box-car kernels, it is  known that, when the number  of channels in the multichannel deconvolution model (\ref{sampl}) is finite, the precision of reconstruction of the response function increases as the number of channels $M$ grow (even when the total number of observations $n$ for all channels $M$ remains constant) and this requires the channels to form a Badly Approximable (BA)  $M$-tuple (see De Canditiis and Pensky (2004, 2007)).  Under the same assumption for the errors, Pensky and Sapatinas~(2009, 2010) showed that the construction of a BA  $M$-tuple for the channels is not needed and a uniform sampling strategy for the channels with the number of channels increasing at a polynomial rate (i.e., $u_l=l/M$, $l=1,2,\ldots,M$, for $M =M_n \geq  (32 \pi/3)n^{1/3}$) suffices to construct an adaptive wavelet estimator that is asymptotically optimal (in the minimax sense),
or near-optimal within a logarithmic factor, in a wide range of Besov balls, when the blurring function $g(\cdot,\cdot)$ is of box-car like kernel (including both the standard box-car kernel and the kernel that appears the estimation of the initial speed of a wave on a finite interval). Example \ref{ex2} showed that a similar result is still possible under long-range dependence with (equispaced channels) $u_l=l/M$, $l=1,2,\ldots,M$, $n^{\te_1} \leq M=M_n \leq n^{\te_2}$, for some $0 \leq \te_1 \leq \te_2 <1$ when $d_l= a_1 u_l + a_2$, $l=1,2,\ldots,M$, $0 \leq a_2 <  1/2$, $0 \leq a_1 + a_2 <  1/2$.

However, in real-life situations, the number of channels $M = M_n$ usually refers to the number of physical devices
and, consequently,  may grow to infinity only at a slow  rate as $n \rightarrow \infty$. When  $M=M_n$ grows slowly as  $n$ increases, (i.e., $M=M_n=o((\ln n)^{\alpha})$ for some $\alpha \geq 1/2$), in the multichannel deconvolution model with  independent and identically distributed Gaussian errors, Pensky and Sapatinas (2011) developed a procedure for the construction of a BA $M$-tuple on a specified interval, of a non-asymptotic length, together with a lower bound associated with this $M$-tuple, which explicitly shows its dependence on $M$ as $M$ is growing. This result was further used for the derivation of upper bounds for the $L^2$-risk of the suggested adaptive wavelet thresholding estimator of the unknown response function and, furthermore, for the choice of the optimal number of channels $M$ which minimizes the $L^2$-risk.
It would be of interest to see whether or not similar upper bounds are possible under long-range dependence. Another avenue of possible research is to consider an analogous minimax study for the functional deconvolution model (i.e., the continuous version of the multichannel deconvolution model (\ref{sampl})) under long range-dependence (e.g., modeling the errors as a fractional Brownian motion) and examine the effect of the convergence rates between the two models, similar to the convergence rate  study of Pensky and Sapatinas (2010) when the errors were considered to be independent and identically distributed Gaussian random variables.


\input{MultichannelProofs}

\section*{Acknowledgements}
Marianna Pensky was supported in part by National Science Foundation
(NSF), grant  DMS-1106564.


\end{document}

%% file: MultichannelProofs.tex

\section{Proofs}
\label{sec:proofs}
\setcounter{equation}{0}

\subsection{Proofs of the Statements in Section \ref{sec:long-range}}

{\bf Proof of Lemma \ref{lemma:assump1}. }
We prove the upper bound only since the  proof of the
lower bound is similar. By
(\ref{eq:spectral-radius})-(\ref{eq:spectral-repr}), and the
definitions of ${\cal H}_N$ and ${\cal F}_d$,
$$
\lambda_{\max} (\bSigma) \le C_{\max} \sup_{h\in {\cal
H}_N}\int_{-\pi}^{\pi}h(\lambda) |\lambda|^{-2d} d\lambda =2C_{\max}
\sup_{h\in {\cal H}_N}\int_{0}^{\pi}h(\lambda) |\lambda|^{-2d}
d\lambda.
$$
Now, we split $\int_{0}^{\pi}=\int_{0}^{\pi/N}+\int_{\pi/N}^{\pi}$.
Since $d <1/2$, for the first integral, we have
$$
\int_0^{\pi/N}h(\lambda) |\lambda|^{-2d} d\lambda \le
N\int_0^{\pi/N} \lambda^{-2d} d\lambda = N\frac{1}{1-2d}
\left(\frac{\pi}{N}\right)^{-2d+1} =\frac{1}{1-2d}N^{2d}.
$$
For the second integral, since $d\ge 0$, we have
$$
\int_{\pi/N}^{\pi}h(\lambda) |\lambda|^{-2d} d\lambda \le
\left(\frac{\pi}{N}\right)^{-2d}\int_{\pi/N}^{\pi}h(\lambda)
 d\lambda \le  \left(\frac{\pi}{N}\right)^{-2d} \int_{0}^{\pi}h(\lambda)
 d\lambda  \le \pi (2\pi)^{-2d} N^{2d}.
$$
This completes the proof of the lemma. \hfill $\Box$

\subsection{Proof of the Minimax Lower Bounds for the $L^2$-Risk}

In order to prove Theorem \ref{th:lower}, we consider two cases: the  dense case and the sparse case,
when the hardest functions to estimate  are, respectively,  uniformly spread over the unit interval $T$
and are represented by only one term in a wavelet expansion. \\

The proof of Theorem \ref{th:lower} is based on Lemma A.1 of Bunea, Tsybakov and Wegkamp~(2007),
which we reformulate here for the case of the $L^2$-risk.

\begin{lemma} \label{Tsybakov} [Bunea, Tsybakov, Wegkamp~(2007), Lemma A.1]
Let $\Omega$ be a set of functions of cardinality $\card(\Omega) \geq 2$, such that \\
(i) $\|f-g\|^2 \geq 4 \delta^2$  for $f, g \in \Omega$, \ $f \neq g$,\\
(ii) the Kullback divergences $K( P_f, P_g)$ between the measures $P_f$ and $P_g$ satisfy the inequality
 $K( P_f, P_g) \leq   \log(\card(\Omega))/16$ for $f, g \in \Omega$.\\
Then, for some absolute constant $C>0$, one has
 \begin{eqnarray*}
 \inf_{T_n} \sup_{f \in \ \Omega} \EE_f  \|T_n -f\|^2  \geq C \delta^2.
 \end{eqnarray*}
\end{lemma}

\noindent
\underline{\bf  The dense case}.
Let  $\bomega$ be the vector with components $\omega_{ k}=\{ 0,1\}$.
Denote the set of all possible vectors  $\bomega$ by  $\Omega= \{ (0, 1)^{2^j}\}$.
Note that the vector  $\bomega$ has $\aleph=2^{j}$ entries  and, hence,   $\card(\Omega)=2^{\aleph}$.
Let $H( \tilde{\bomega}, \bomega)=  \sum^{2^j-1}_{k=0} \II \left(\tilde{\bomega}_{k} \neq \omega_{k}\right)$
be the Hamming distance between the binary sequences $ \bomega$ and $ \tilde{\bomega}$.
Then, the Varshamov-Gilbert Lemma (see, e.g., Tsybakov (2008), p. 104) states that one can choose
a subset $\Omega_1$ of $\Omega$, of cardinality at least $2^{\aleph/8}$, such that $H(\tilde{\bomega}, \bomega) \geq \aleph/8$
for any $\bomega , \tilde{\bomega} \in \Omega_1$.\\

Consider two arbitrary sequences  $\bomega , \tilde{\bomega} \in \Omega_1$ and  the
functions $f_{j}$ and $\tilde{f}_{j}$  given by
\begin{eqnarray*}
f_{j}(t)= \rho_{j}  \sum^{2^{j}-1}_{k=0} \omega_{k}\psi_{jk}(t) \quad \text{and} \quad
\tilde{f}_{j}(t)= \rho_{j}  \sum^{2^{j}-1}_{k=0}\tilde{ \omega}_{k}\psi_{jk}(t), \quad t \in T.
\end{eqnarray*}
Choose $\rho_{j} = A  2^{- j(s+1/2)}$, so that $f_{j}, \tilde{f}_{j} \in  B^{s}_{p, q}(A)$.
Then, calculating the $L^2$-norm difference of $f_{j}$ and $\tilde{f}_{j}$, we obtain
\begin{eqnarray*}
 \| \tilde{f}_{j}-f_{j}\|^2 =  \rho^2_{j}\|  \sum^{2^{j}-1}_{k=0}
(\tilde{\omega}_{k}- \omega_{k})\psi_{jk} \|^2
= \rho^2_{j} H(\tilde{ \bomega}, \bomega) \geq 2^j \rho^2_{j}/8.
\end{eqnarray*}
Hence, we get $4 \delta^2  = 2^j \rho^2_{j}/8$ in  condition $(i)$ of Lemma \ref{Tsybakov}. \\

In order to apply Lemma \ref{Tsybakov}, one needs to also verify condition $(ii)$.
For $f_{\bomega}$ with $\bomega \in \Omega$, denote by $\bh_{l, \bomega}$ and  $\bh_{l, \tilde{\bomega}}$,  the vectors
with components, respectively,
\begin{eqnarray*}
h_{\bomega}(u_l, t_i) &=&    g(u_l, t_i - \cdot) * f_{\bomega}(\cdot), \quad i= 1,2,\ldots,N, \\
h_{\tilde{\bomega}}(u_l, t_i) &=&    g(u_l, t_i - \cdot) * f_{\tilde{\bomega}}(\cdot), \quad  i= 1,2,\ldots,N.
 \end{eqnarray*}
Then,
\begin{eqnarray*}
K(P_{\bomega}, P_{\tilde{\bomega}}) & = &
0.5 \sum^M_{l=1}    (\bh_{l, \bomega}  - \bh _{l, \tilde{\bomega}})^T (\bSigma^{(l)})^{-1} (\bh_{l, \bomega}  - \bh _{l, \tilde{\bomega}}) \\
& \leq &
0.5 \sum^M_{l=1}  \lambda_{\max} ( (\bSigma^{(l)})^{-1})  \|\bh_{l, \bomega}  - \bh _{l, \tilde{\bomega}}\|^2.
\end{eqnarray*}
Now, since $\bomega$ and $\tilde{\bomega}$ are binary vectors, using Plancherel's formula, \fr{delta1}, and
the fact that $| \psi_{jk, m}|\leq 2^{-j/2}$, we derive that, under Assumptions A1 and A2,
\begin{eqnarray*}
K(P_{\bomega}, P_{\tilde{\bomega}})
& \leq &
0.5 NM \, \rho^2_j\   \sum_{m \in C_j} \frac{1}{M}  \sum^M_{l=1}|g_m(u_l)|^2 K_1^{-1} N^{-2 d_l} \\
& \leq &
2 \pi K_1^{-1} n  2^{j} \rho^2_j\ \Delta_1(j,n)  \leq 2 \pi A^2 K_1^{-1} n 2^{-2js}\  \Delta_1(j,n),\\
\end{eqnarray*}
where $\Delta_1(j,n)$ is defined by \fr{Delt}. \\

Direct calculations yield that, under Assumptions A1, A2  and \fr{ns_prop1},
for some constants $c_3 >0$ and $c_4 >0$, independent of $n$,
\be
\Delta_1(j,n) \leq \lfi
\begin{array}{ll}
c_3\ \eps_n^{-1}\ 2^{2 \nu j} j^{\lam_2}, &
\mbox{if}\;\;\; \alpha_1=\alpha_2=0, \\
& \label{delta1}\\
c_4\ \eps_n^{-1}\ 2^{2 \nu_1 j} j^{\lam_2}\ \exp \lfi \alpha_1 \lkr \frac{8\pi}{3} \rkr^\beta
2^{j \beta} \rfi,  & \mbox{if}\;\;\;  \alpha_1 \alpha >0. 
\end{array} \right.
\ee
Apply  now Lemma \ref{Tsybakov} with  $j$ such that
$$
2 \pi A^2 K_1^{-1} n 2^{-2js}\  \Delta_1(j,n) \leq 2^j\, \ln 2/16,
$$
i.e.,
$$
2^j   \asymp \left\{ \begin{array}{ll}
\left[  n^* (\ln n^*)^{-\lambda_2} \right]^{\frac{1}{2s +2 \nu +1}}, & \mbox{if}\ \ \beta=0,\\
(\ln n^* )^{1/\beta}, & \mbox{if}\ \ \beta >0,
 \end{array} \right.
$$
to obtain
\be \label{lbdd}
\delta^2=  \left\{ \begin{array}{ll}
 \left[  n^* (\ln n^*)^{-\lambda_2} \right]^{-\frac{2s}{2s+2\nu +1}}, & \mbox{if}\ \ \ \beta = 0, \\
  (\ln n^* )^{-2s/\beta}, & \mbox{if} \ \ \beta > 0.
 \end{array} \right.
\ee


\noindent
\underline{\bf  The sparse case}.
Let the functions $f_j$ be of the form $f_j(t)= \rho_j  \psi_{jk}(t)$, $t \in T$, and denote
$$\Omega= \{ f_j(t)= \rho_j  \psi_{jk}(t): \; k= 0,1,\ldots, 2^j-1, \; f_0= 0 \}.$$
Thus, $\card(\Omega)=2^{j}$.
Choose now $\rho_j = A 2^{-js'}$, so  that $f_{j} \in B^{s}_{p, q}(A)$.
It is easy to check that, in this case, one has $4\delta^2=  \rho^2_{j}$ in  Lemma \ref{Tsybakov},
and that
\begin{eqnarray*}
K(P_{\omega}, P_{\tilde{\omega}}) & \leq &
2 \pi A^2 K_1^{-1} n 2^{-2js'}\  \Delta_1(j,n).
\end{eqnarray*}
With
$$
2^j   \asymp \left\{ \begin{array}{ll}
\left[  n^* (\ln n^*)^{-\lambda_2-1}  \right]^{\frac{1}{2s' +2 \nu}}, & \mbox{if}\ \ \beta=0,\\
(\ln n^* )^{1/\beta}, & \mbox{if}\ \ \beta >0,
 \end{array} \right.
$$
we then obtain that $K(P_{\omega}, P_{\tilde{\omega}}) \leq 2 \pi A^2 K_1^{-1} n 2^{-2js'}\  \Delta_1(j,n)$
and
\be \label{lbds}
\delta^2=  \left\{ \begin{array}{ll}
  \left[ \frac{n^*}{ (\ln n^*)^{\lambda_2+1}}\right]^{-\frac{2s'}{2s'+2\nu}}, & \mbox{if}\ \ \beta=0, \\
  (\ln n^* )^{-2s'/\beta}, & \mbox{if} \ \beta >0.
 \end{array} \right.
\ee

Recall that $s^*= \min\{ s, s'\}$. By noting that
\be
2s/(2s +2\nu +1) \leq 2s^*/(2s^*+ 2\nu),   \quad  \mbox{if} \quad   \nu(2-p) \leq ps^*,
\ee
we then choose the highest of the lower bounds in \fr{lbdd} and \fr{lbds}.
This completes the proof of the theorem. \hfill $\Box$


\subsection{Proof of the Minimax Upper Bounds for the $L^2$-Risk.}

 We start with proofs of  Lemmas \ref{l:coef} and \ref{l:deviation}.
\\

{\bf Proof of Lemma \ref{l:coef}.}
First, consider model \fr{sampl}. Then, using \fr{hgf2}, \fr{fmh},  \fr{alkandblk} and \fr{coefest}, one has
\begin{eqnarray*}
\widehat{a}_{j_0k} - a_{j_0k}= \sum_{m \in C_{j_0}}\left( \widehat{f}_m - f_m \right) \overline{\varphi_{mj_0k}}, \quad
\widehat{b}_{jk}-b_{jk}=  \sum_{m \in C_{j}}\left( \widehat{f}_m- f_m \right) \overline{\psi_{mjk}},
\end{eqnarray*}
where
\be
\widehat{f}_m - f_m= \frac{1}{\sqrt{N}} \left( \sum^M_{l=1}N^{-2 d_l} \overline{g_m(u_l)}z_{lm}  \right)/
\left( \sum^{M}_{l=1}N^{-2d_l} |g_m(u_l)|^2 \right).
\ee
Consider vector $\bV^{(l)}$ with components
$$V^{(l)}_m = N^{- 2 d_l} \psi_{mjk}   g_m(u_l)  \lkv
\sum^{M}_{l=1} N^{-2d_l}|g_{m}(u_l)|^2 \rkv^{-1}.$$ It is easy to see
that, due to $|\psi_{mjk}| \leq 2^{-j/2}$ and the definition of $C_j$,
\begin{eqnarray*}
\|\bV^{(l)}\|^2 & = & N^{- 4 d_l} \ \sum_{m  \in C_{j}} | \psi_{mjk}|^2
|g_m(u_l)|^2 \lkv \sum^{M}_{l=1} N^{-2d_l}|g_{m}(u_l)|^2 \rkv^{-2}\\
& \leq   & 4 \pi |C_j|^{-1}\ N^{- 4 d_l} \sum_{m \in C_{j}}
|g_m(u_l)|^2 \lkv \sum^{M}_{l=1} N^{-2d_l}|g_{m}(u_l)|^2 \rkv^{-2}.
\end{eqnarray*}
Define
$$
v_m=\sum^{M}_{l=1} N^{-2d_l}|g_{m}(u_l)|^2=M\tau_1(m,n)\;.
$$
Hence,
$$
\|\bV^{(l)}\|^2 \leq 4 \pi |C_j|^{-1}\ N^{- 2 d_l}  N^{- 2
d_l}\sum_{m \in C_{j}} |g_m(u_l)|^2 v_m^{-2}.
$$

Using Assumption A1, since  $z_{lm}$ are
independent for different $l'$s,  we obtain
\begin{eqnarray*}
\EE|\widehat{b}_{jk}-b_{jk}|^2 & =  & \frac{1}{N} \sum_{m_1, m_2 \in
C_{j}}\overline{\psi}_{m_1jk} \psi_{m_2jk}
\sum^M_{l=1} N^{-4d_l} v_{m_1}^{-1}v_{m_2}^{-1}  \overline{g_{m_1}(u_l)}
g_{m_2}(u_l)
\Cov  \left(z_{lm_1}, \overline{z}_{lm_2}  \right)\\
& = &
\frac{1}{N} \sum_{l=1}^M  \  \overline{\bV^{(l)}}^T \bSigma^{(l)}
\bV^{(l)} \\
& \leq   &
\frac{1}{N} \sum_{l=1}^M  \  \lambda_{\max} (\bSigma^{(l)}) \|
\bV^{(l)}\|^2 \\
&\leq & 4\pi K_2 \, |C_j|^{-1} N^{-1} \sum_{l=1}^M N^{- 2
d_l}\sum_{m \in C_{j}} |g_m(u_l)|^2 v_m^{-2}\\
&= & 4\pi K_2 \, |C_j|^{-1} N^{-1}\sum_{m \in C_{j}} v_m^{-2}
\sum_{l=1}^M N^{- 2 d_l} |g_m(u_l)|^2 = 4\pi K_2 \, |C_j|^{-1}
N^{-1}\sum_{m \in C_{j}} v_m^{-1}  ,
\end{eqnarray*}
so that
\begin{eqnarray*}
\EE|\widehat{b}_{jk}-b_{jk}|^2  & \leq   &
C n^{-1}  |C_j|^{-1}  \sum_{m \in C_{j}}\left[ \tau_1 (m,n) \right]^{-1}.
\end{eqnarray*}
(One can obtain an  upper bound  for $\EE| \widehat{a}_{j_0k} -
a_{j_0k} |^2$ by following similar arguments.)


In order  to prove \fr{hbb},  define
$$
B_l=  N^{- 2d_l} \lkv \sum^{M}_{l=1}N^{-2d_l} |g_m(u_l)|^2 \rkv^{-1}.
$$
Note that
\begin{eqnarray*}
\EE \left( z_{lm_{1}}z_{lm_2}z_{lm_3}z_{lm_4} \right)
& \leq &  \left[ \Pi^{4}_{i=1} \EE | z_{m_il} |^4 \right]^{1/4}.
\end{eqnarray*}
Consequently,  using Assumption A1, the fact that  $z_{lm}$ are independent for different $l'$s,
and that  $\EE | z_{ml}|^4=3 \left[ \EE| z_{ml} |^2 \right]^2$ for standard (complex-valued) Gaussian random variables $z_{ml}$, one obtains
\begin{eqnarray*}
 \EE|\widehat{b}_{jk}-b_{jk}|^4
& = &
 O \left(  N^{-2}\   \sum^M_{l=1}  B_l^4  \left[ \sum_{m \in C_{j}} |\psi_{mjk}| |g_{m_2}(u_l)| \left(\EE | z_{ml} |^4 \right)^{1/4} \right]^4  \right) \\
 & + &
O \left(\left[N^{-1}\  \sum^M_{l=1}  B_l^2  \sum_{m_1, m_2 \in C_{j}} \overline{\psi}_{m_1jk} \psi_{m_2jk}
\overline{g_{m_1}(u_l)}  g_{m}(u_l) \Cov  \left(z_{lm_1}, \overline{z}_{lm_2}  \right)  \right]^2 \right) \\
  & = &
O \left( N^{-2}\   {\sum^M_{l=1} B_l^4 \left[\sum_{m \in C_{j}} | \psi_{mjk}|^2 |g_{m}(u_l)|^2 \sum_{m \in C_{j}}  \EE | z_{ml} |^2\right]^2} \right) \\
  & + &
 O \left( \left[  n^{-1} |C_j|^{-1}\  \sum_{m \in C_{j}} \left[ \tau_1 (m,n) \right]^{-1} \right]^2 \right)
\end{eqnarray*}
Since
$\sum_{m \in C_{j}}  \EE | z_{ml} |^2 = O(|C_{j}|)$, one derives
\begin{eqnarray}
 \EE|\widehat{b}_{jk}-b_{jk}|^4
& = &
 O \left(|C_j|^{-1}    \sum_{m \in C_{j}}\left[\frac{1}{M^3} \frac{ \tau_2 (m,n)}{[\tau_1 (m,n)]^4}   \right] +
\frac{\Delta^2_1 (j,n)}{n^2}  \right) \nonumber  \\
& = & O \left(  M^{-3}\Delta_2 (j,n) + n^{-2}\Delta^2_1 (j,n)  \right).  \label{formo}
\end{eqnarray}
It is straightforward to show that, when $\alpha_1=\alpha_2=0$, one has
$$\Delta_2 (j,n)= O \left(  2^{6j\nu} j^{3\lambda_1}\varepsilon^{-3}_n \right).$$
Thus, using   \fr{delta1} and   the fact that $2^j \leq 2^{J-1} < (n^*)^{1/(2\nu +1)}$, \fr{formo} can be rewritten as
\begin{eqnarray*}
 \EE|\widehat{b}_{jk}-b_{jk}|^4
  & = &
O \left(  2^{6\nu j} j^{3\lambda_1}\varepsilon^{-3}_n M^{-3}  + 2^{4j\nu} j^{2\lam_1}\varepsilon^{-2}_n n^{-2}  \right) \\
& = &  O \left( n^3 \left( \ln n \right)^{3\lam_1} \left( n^*\right)^{-3/(2\nu +1)} \right).
\end{eqnarray*}
Hence, \fr{hbb} follows. This completes the proof of the lemma. \hfill $\Box$
\\


{\bf Proof of Lemma \ref{l:deviation}. } Consider a set of vectors
\begin{eqnarray*}
 \Omega_{jr}= \left \{ v_k, k \in U_{jr}:  \sumku | v_k |^2 \leq 1   \right \}
   \end{eqnarray*}
   and a centered Gaussian process
$$
Z_{jr}=  \sumku v_k ( \widehat{b}_{jk} - b_{jk} ).
$$
   Note that
   $$\sup_v Z_{jr}(v)= \sqrt{ \sumku |  \widehat{b}_{jk} - b_{jk} |^2 }.$$

   We shall apply below a
lemma of Cirelson, Ibragimov and Sudakov~(1976) which states that, for any $x >0$,
   \be \label{CarlIbSud}
   \Pr \left( \sumku |  \widehat{b}_{jk} - b_{jk} |^2 \geq (x + B_1) \right) \leq \exp \left( -\frac{x^2}{2B_2} \right),
   \ee
  where,
   \begin{eqnarray*}
   B_1   & = & \EE \left[ \sqrt{ \sumku |  \widehat{b}_{jk} - b_{jk} |^2 } \right]
     \leq  \frac{\sqrt{c_1} 2^{j\nu} j^{\lambda_1/2} \sqrt{\ln n}}{\sqrt{ n^*}}
   \end{eqnarray*}
   with $c_1$ defined in \fr{delta1},  and
   \begin{eqnarray*}
   B_2 = \sup_{v \in \Omega_{jr}} \Var (Z_{jr}(v))= \sup_{v \in \Omega_{jr}} \EE  |\sumku v_k ( \widehat{b}_{jk} - b_{jk} )     |^2.
      \end{eqnarray*}

Denote
$$
w_{jm}= \sumku v_k  \psi_{mjk} \lkv \sum^{M}_{l=1} N^{-2d_l}|g_{m}(u_l)|^2 \rkv^{-1},\quad m \in C_j.
$$
Then, under Assumption A2  with $\alpha_1 = \alpha_2 = 0$,
using argument similar to the proof of \fr{ha}, one obtains
    \begin{eqnarray*}
       B_2 & = &
\sup_{v \in \Omega_{jr}}    \left \{N^{-1}\    \sum_{m_1, m_2 \in C_{j}} \overline{w_{jm_1}} w_{jm_2}\
\EE \left[  \sum^M_{l=1} N^{-4d_l} \overline{g_{m_1}(u_l)} g_{m_2}(u_l) z_{lm_1} \overline{z}_{lm_2}   \right]  \right\}\\
  & \leq &
\sup_{v \in \Omega_{jr}}   N^{-1}\  \sum_{l=1}^M  N^{-4 d_l} \lam_{\max}( \bSigma^{(l)})
\sum_{m \in C_{j}} |w_{jm}  g_m (u_l)|^2 \\
 & \leq &
K_3 n^{-1}  \sup_{v \in \Omega_{jr}}    \left \{  \sum_{m \in C_{j}} |w_{jm}  |^2 \left[ \tau_1 (m,n) \right]^{-1}  \right \}
\leq 4 \pi C^*_3   2^{2j\nu} j^{\lambda_1}\,(n^*)^{-1},
\end{eqnarray*}
where $C^*_3 = (K_3)^{-1} (\ln 2)^{\lambda_1} (2\pi /3)^{2\nu}$.\\

Apply now inequality \fr{CarlIbSud} with  $x$ such that   $ x^2=2 B_2  \kappa  \ln n$, and note that
   \begin{eqnarray*}
  (x + B_1)^2 =  (n^*)^{-1}  2^{2 j \nu} j^{\lambda_1}  \ln n \
\lkr \sqrt{c_1} +  \sqrt{   8 \pi \kappa K_3^{-1}  (\ln 2)^{\lambda_1}(2\pi/3)^{2\nu} } \rkr^2
     \end{eqnarray*}
and
\begin{eqnarray*}
     \mu^2  \geq 4 (1- h_1 )^{-1} \left ( \sqrt{c_1} + \sqrt{   8 \pi \kappa K_3^{-1}  (\ln 2)^{\lambda_1}(2\pi/3)^{2\nu} } \right)^2,
\end{eqnarray*}
which guarantees  \fr{CarlIbSud}. This completes the proof of the lemma. \hfill $\Box$
\\


{\bf Proof of Theorem \ref{th:upper}.}
Direct calculations yield that under Assumptions A1, A2  and \fr{ns_prop2},
for some constants $c_1 >0$ and $c_2 >0$, independent of $n$, one has
\be
\Delta_1(j,n) \leq \lfi
\begin{array}{ll}
c_1\ \eps_n^{-1}\ 2^{2 \nu j} j^{\lam_1}, &
\mbox{if}\;\;\; \alpha_1=\alpha_2=0, \\
& \label{delta1upper}\\
c_2\ \eps_n^{-1}\ 2^{2 \nu_1 j} j^{\lam_1}\ \exp \lfi \alpha_1 \lkr \frac{8\pi}{3} \rkr^\beta
2^{j \beta} \rfi,  & \mbox{if}\;\;\;  \alpha_1 \alpha >0. 
\end{array} \right.
\ee
Using \ref{delta1upper}, the proof of this theorem is now almost identical to the proof of Theorem 2 in Pensky and Sapatinas~(2010).
This completes the proof of the theorem. \hfill $\Box$


\subsection{Proofs of the Statement in Section \ref{sec:examples}. }

{\bf Proof of Lemma \ref{lem:boxcar}}.
Below we consider only the case of $a_1 >0$. Validity of the satement for $a_1 =0$ follows from Pensky and Sapatinas (2010).

By direct calculations, one obtains that
$$
\tau_1 (m,n) = M^{-1} (4 \pi^2 m^2)^{-1} N^{-2a_2} \sum_{l=1}^M  q^2(l/M)\,  \sin^2(2 \pi m l M^{-1}) N^{- 2 a_1 l/M}.
$$
Therefore,
\be \label{eq:lowup}
 (4 \pi^2 m^2)^{-1} q_1^2 \  N^{-2a_2}  S(m,n)  \leq \tau_1 (m,n) \leq  (4 \pi^2 m^2)^{-1} q_2^2 \  N^{-2a_2}  S(m,n),
\ee
where
$$
S(m,n) = M^{-1}   \sum_{l=1}^M  \sin^2(2 \pi m l M^{-1}) N^{- 2 a_1 l/M}.
$$
Denote $p=N^{-2a_1/M}$, $x=4 \pi m M^{-1}$ and note that, as $n \to \infty$,
$$p^M = N^{-2 a_1} \rightarrow 0$$
and
\begin{eqnarray} \label{p_expr}
p &=& \exp \lkr -2 a_1 M^{-1} \ln N  \rkr \nonumber \\
&=& 1 -  2 a_1 M^{-1} \ln N +
2   a_1^2 M^{-2}  \ln^2  N  + o(M^{-2}  \ln^2  N),
\end{eqnarray}
since $M^{-1} \ln N \to 0$ as $n \to \infty$.\\

Using the fact that $\sin^2(x/2) = (1- \cos x)/2$ and formula 1.353.3 of Gradshtein \& Ryzhik (1980), we obtain
$$
S(m,n) = \frac{1}{M} \lkv \frac{1 -p^M}{1-p} - \frac{1 - p\cos x - p^M \cos(Mx) + p^{M+1}\cos((M-1)x)}
{1 - 2p \cos x + p^2} \rkv.
$$
Since $m$ is an integer and $x=4 \pi m M^{-1}$,
$$\cos(Mx) =1, \quad \sin(Mx)=0, \quad \cos((M-1)x) = \cos x.$$
Therefore, simple algebraic transformations yield
$$
S(m,n) =  \frac{p (p+1)(1 -p^M)(1 - \cos x)}{M(1-p)[(1-p)^2 + 2p(1- \cos x)]}
$$
The asymptotic expansion \fr{p_expr} for $p$ as $n \to \infty$, leads to
\be \label{asymp1}
\frac{ (1 -p^M)}{M(1-p)}  \approx \frac{1 - N^{-2 a_1}}{4 a_1 \ln N (1 - a_1 M^{-1} \ln N)},
\ee
so that, if $N$ is large enough, due to $p<1$, one obtains an upper bound for $S(m,n)$:
\be \label{Smn_upper}
S(m,n) = \frac{ (1 -p^M)}{M(1-p)} \lkv \frac{(1-p)^2}{p (p+1)(1 - \cos x)} + \frac{2}{p+1} \rkv^{-1}
\leq \frac{1}{2 a_1 \ln N}.
\ee

In order to obtain  a lower bound  for $S(m,n)$, we note that  for $N$ large enough,
one has $1/2 < p < 1$. Consider the following two cases: $x \geq \pi/3$ and $x < \pi/3$.
If $x \geq \pi/3$, then $\cos x \leq 1/2$ and
$$
F(p,x) = \frac{(1-p)^2}{p(p+1)(1 - \cos x)} + \frac{2}{p+1} \leq 2,
$$
If $x < \pi/3$, we can use the fact that
$1 - \cos x    = 2\sin^2 (x/2) \geq  3 x^2/8,$
so that
$$
F(p,x) \leq \frac{4}{3} \lkv 1 +  \frac{8 (1-p)^2}{3 x^2} \rkv \leq
\frac{4}{3} \lkv 1 + \frac{2 a_1^2 \ln^2 N}{3 \pi^2 m^2} \rkv
$$
for $N$ large enough. \\

Since  $|m| = C_m 2^j  > C_m C_0  \ln n$ for some $\delta>0$ and
$\ln n \geq (1 - \theta_1)^{-1} \ln N$ due to assumption \fr{eq:M_cond}, one has
$m^2 \geq C_m C_0 (1 - \te_1)^{-1}  \ln^2 N$ and
 \be \label{Smn_lower}
S(m,n)   \geq C (\ln N)^{-1}.
\ee
Observe now that $\ln N \asymp \ln n$.
This completes the proof of the theorem. \hfill $\Box$